\documentclass[11pt]{article}%
\usepackage{amsmath}
\usepackage{amsfonts}
\usepackage{amssymb}
\usepackage{graphicx}
\usepackage{a4wide}
\newtheorem{theorem}{Theorem}

\newtheorem{corollary}[theorem]{Corollary}

\newtheorem{example}[theorem]{Example}

\newtheorem{lemma}[theorem]{Lemma}

\newtheorem{proposition}[theorem]{Proposition}
\newtheorem{remark}[theorem]{Remark}

\begin{document}

\title{On the global shape of convex functions on locally convex spaces}
\author{C. Z\u{a}linescu\\Octav Mayer Institute of Mathematics\\Iasi Branch of Romanian Academy, Iasi, Romania}
\maketitle

\begin{abstract}
In the recent paper \cite{Aza:19} D Azagra studies the global shape of
continuous convex functions defined on a Banach space $X$. More precisely,
when $X$ is separable, it is shown that for every continuous convex function
$f:X\rightarrow\mathbb{R}$ there exist a unique closed linear subspace $Y$ of
$X$, a continuous function $h:X/Y\rightarrow\mathbb{R}$ with the property that
$\lim_{t\rightarrow\infty}h(u+tv)=\infty$ for all $u,v\in X/Y$, $v\neq0$, and
$x^{\ast}\in X^{\ast}$ such that $f=h\circ\pi+x^{\ast}$, where $\pi
:X\rightarrow X/Y$ is the natural projection. Our aim is to characterize those
proper lower semi\-continuous convex functions defined on a locally convex
space which have the above representation. In particular, we show that the
continuity of the function $f$ and the completeness of $X$ can be removed from
the hypothesis of Azagra's theorem.

\end{abstract}

\section{Preliminary notions and results}

In the sequel $X$ is a nontrivial real separated locally convex space (lcs for
short) with topological dual $X^{\ast}$ endowed with its weak$^{\ast}$
topology (if not explicitly mentioned otherwise); for $x\in X$ and $x^{\ast
}\in X^{\ast}$ we set $\left\langle x,x^{\ast}\right\rangle :=x^{\ast}(x)$. In
some statements $X$ will be a real normed vector space (nvs for short), or
even a Hilbert space, in which case $X^{\ast}$ will be identified with $X$ by
Riesz theorem. For $E$ a topological vector space and $A\subset E$, we denote
by $\overline{A}$ (or $\operatorname*{cl}A$) and $\operatorname*{span}A$ the
closure and the linear hull of $A$, respectively; moreover, $\overline
{\operatorname*{span}}A:=\overline{\operatorname*{span}A}$. In particular,
these notations apply for the subsets of $X^{\ast}$ which is endowed with the
weak-star topology by default; when $X$ is a normed vector space, the
norm-closure of $B\subset X^{\ast}$ is denoted by $\operatorname*{cl}%
_{\left\Vert \cdot\right\Vert }B$.

The \emph{domain} of the function $f:X\rightarrow\overline{\mathbb{R}%
}:=\mathbb{R}\cup\{-\infty,\infty\}$ is the set $\operatorname*{dom}f:=\{x\in
X\mid f(x)<\infty\}$. The function $f$ is \emph{proper} if
$\operatorname*{dom}f\neq\emptyset$ and $f(x)>-\infty$ for all $x\in X$; $f$
is \emph{convex} if $\operatorname*{epi}f:=\{(x,t)\in X\times\mathbb{R}\mid
f(x)\leq t\}$ is convex. Of course, $f$ is \emph{lower semi\-continuous} (lsc
for short) iff $\operatorname*{epi}f$ is a closed subset of $X\times
\mathbb{R}$. By $\Gamma(X)$ we denote the class of proper lower
semi\-continuous convex functions $f:X\rightarrow\overline{\mathbb{R}}$.
Having $f:X\rightarrow\overline{\mathbb{R}}$, its \emph{conjugate} function is%
\[
f^{\ast}:X^{\ast}\rightarrow\overline{\mathbb{R}},\quad h^{\ast}(x^{\ast
}):=\sup\left\{  \left\langle x,x^{\ast}\right\rangle -h(x)\mid x\in
X\right\}  \quad(x^{\ast}\in X^{\ast}),
\]
while its \emph{sub\-differential} is the set-valued function $\partial
f:X\rightrightarrows X^{\ast}$ with%
\[
\partial f(x):=\{x^{\ast}\in X^{\ast}\mid\left\langle x^{\prime}-x,x^{\ast
}\right\rangle \leq f(x^{\prime})-f(x)\ \forall x^{\prime}\in X\}
\]
if $f(x)\in\mathbb{R}$ and $\partial f(x):=\emptyset$ otherwise. By
\cite[Th.\ 2.3.3]{Zal:02}, $f^{\ast}\in\Gamma(X^{\ast})$ and $(f^{\ast}%
)^{\ast}=f$ ($X^{\ast}$ being endowed, as mentioned above, with the weak-star
topology $w^{\ast}$) whenever $f\in\Gamma(X)$; in particular
$\operatorname*{dom}f^{\ast}\neq\emptyset$. Moreover, for $f\in\Gamma(X)$ one
has $x^{\ast}\in\partial f(x)$ iff $x\in\partial f^{\ast}(x^{\ast})$ iff
$f(x)+f^{\ast}(x^{\ast})=\left\langle x,x^{\ast}\right\rangle $.

A central notion throughout this note is that of recession function. So,
having $f\in\Gamma(X)$, its \emph{recession function} $f_{\infty}$ is
(equivalently) defined by
\[
f_{\infty}:X\rightarrow\overline{\mathbb{R}},\quad f_{\infty}(u):=\lim
_{t\rightarrow\infty}\frac{f(x_{0}+tu)-f(x_{0})}{t},
\]
where $x_{0}\in\operatorname*{dom}f$ is arbitrary. The function $f_{\infty}$
is a proper lsc sublinear function having the property
\begin{equation}
f(x+u)\leq f(x)+f_{\infty}(u)\quad\forall x\in\operatorname*{dom}f,\ \forall
u\in X \label{r-a0}%
\end{equation}
(see \cite[Eq.\ (2.28)]{Zal:02}); moreover,
\begin{equation}
f_{\infty}(u)=\sup_{x^{\ast}\in\operatorname*{dom}f^{\ast}}\left\langle
u,x^{\ast}\right\rangle ~\forall u\in X\quad\text{and}\quad\partial f_{\infty
}(0)=\overline{\operatorname*{dom}f^{\ast}} \label{r-a5}%
\end{equation}
(see \cite[Exer.\ 2.23 and Th.\ 2.4.14]{Zal:02}). In particular (see also
\cite[Th.\ 2.4.14]{Zal:02}), one has
\begin{equation}
\partial g(0)=\{x^{\ast}\in X^{\ast}\mid x^{\ast}\leq g\},\quad g^{\ast}%
=\iota_{\partial g(0},\quad\text{and}\quad g=g_{\infty}=\sup\nolimits_{x^{\ast
}\in\partial g(0)}x^{\ast}, \label{r-a5b}%
\end{equation}
where $\iota_{A}:V\rightarrow\overline{\mathbb{R}}$ denotes the
\emph{indicator function} of $A\subset V$, being defined by $\iota_{A}(v):=0$
for $v\in A$ and $\iota_{A}(v):=\infty$ for $v\in V\setminus A$. Hence
$\partial g(0)\neq\emptyset$.

Recall that the mapping $0<t\mapsto\frac{f(x_{0}+tu)-f(x_{0})}{t}\in
\overline{\mathbb{R}}$ is nondecreasing for $f:X\rightarrow\overline
{\mathbb{R}}$ a proper convex function, $x_{0}\in\operatorname*{dom}f$ and
$u\in X$. Moreover, for such a function and $x,u\in X$, the mapping
$\varphi_{x,u}:\mathbb{R}\rightarrow\overline{\mathbb{R}}$ with $\varphi
_{x,u}(t):=f(x+tu)$, one of the following alternatives holds: 1)
$\varphi_{x,u}$ is nonincreasing on $\mathbb{R}$, 2) $\varphi_{x,u}$ is
nondecreasing on $\mathbb{R}$, 3)~there exists $t_{0}\in\mathbb{R}$ such that
$\varphi_{x,u}$ is nonincreasing on $]-\infty,t_{0}]$ and nondecreasing on
$[t_{0},\infty\lbrack$; moreover, there exists $\gamma_{x,u}:=\lim
_{t\rightarrow\infty}f(x+tu)\in\overline{\mathbb{R}}$.

\begin{lemma}
\label{lem-Fact a} Let $f\in\Gamma(X)$ and $u\in X\setminus\{0\}$. The
following assertions are equivalent:

\emph{(a)} $\exists x_{0}\in\operatorname*{dom}f$, $\exists M\in\mathbb{R}$,
$\forall t\in\mathbb{R}_{+}:f(x_{0}+tu)\leq M$;

\emph{(b)} $\forall x\in\operatorname*{dom}f$, $\exists M\in\mathbb{R}$,
$\forall t\in\mathbb{R}_{+}:f(x+tu)\leq M$;

\emph{(c)} $f_{\infty}(u)\leq0$.

Consequently, the following assertions are equivalent:

\emph{(a')} $\forall x\in\operatorname*{dom}f:$ $\lim\nolimits_{t\rightarrow
\infty}f(x_{0}+tu)=\infty$;

\emph{(b')} $\exists x_{0}\in\operatorname*{dom}f:\lim\nolimits_{t\rightarrow
\infty}f(x+tu)=\infty$;

\emph{(c')} $f_{\infty}(u)>0$.
\end{lemma}

Proof. (c) $\Rightarrow$ (b) Take $x\in\operatorname*{dom}f$; then, by
(\ref{r-a0}), $f(x+tu)\leq f(x)+f_{\infty}(tu)=f(x)+tf_{\infty}(u)\leq
f(x)=:M$ for $t\geq0$.

(b) $\Rightarrow$ (a) The implication is obvious.

(a) $\Rightarrow$ (c) Since $t^{-1}\left[  f(x_{0}+tu)-f(x_{0})\right]  \leq
t^{-1}\left[  M-f(x_{0})\right]  $ for $t>0$, one has $f_{\infty}(u)\leq
\lim\nolimits_{t\rightarrow\infty}t^{-1}\left[  M-f(x_{0})\right]  =0$.

Observe that $\rceil$(c') coincides with (c), $\rceil$(b') is equivalent to
(b), and $\rceil$(a') is equivalent to (a). Hence, from the first part, we get
(a')$~\Leftrightarrow$ (b')$~\Leftrightarrow$ (c'). \hfill$\square$

\medskip Having in view the statements of Theorems 5 and 6 in \cite{Aza:19},
it is worth observing that for $x_{0}\in\operatorname*{dom}f$, $u\in X$ and
$u^{\ast}\in X^{\ast}$ one has
\begin{equation}
f_{\infty}(\pm u)=\left\langle \pm u,u^{\ast}\right\rangle \Longleftrightarrow
\left[  f(x_{0}+tu)-f(x_{0})-\left\langle tu,u^{\ast}\right\rangle
=0\ \ \forall t\in\mathbb{R}\right]  . \label{r-a4}%
\end{equation}

\noindent Indeed, the implication \textquotedblleft$\Leftarrow$" is obvious.
Assume that $f_{\infty}(\pm u)=\left\langle \pm u,u^{\ast}\right\rangle $
($\Leftrightarrow f_{\infty}(tu)=\left\langle tu,u^{\ast}\right\rangle $ for
all $t\in\mathbb{R}$). Using (\ref{r-a0}) we get
\[
f(x_{0}+tu)\leq f(x_{0})+f_{\infty}(tu)=f(x_{0})+\left\langle tu,u^{\ast
}\right\rangle ,\text{ \ }f(x_{0}-tu)\leq f(x_{0})-\left\langle tu,u^{\ast
}\right\rangle \quad\forall t\in\mathbb{R}.
\]
Since $x_{0}=\tfrac{1}{2}(x_{0}+tu)+\tfrac{1}{2}(x_{0}-tu)$, from the
convexity of $f$ and the previous inequalities we get%
\[
f(x_{0})\leq\tfrac{1}{2}f(x_{0}+tu)+\tfrac{1}{2}f(x_{0}-tu)\leq\tfrac{1}%
{2}[f(x_{0})+\left\langle tu,u^{\ast}\right\rangle ]+\tfrac{1}{2}%
[f(x_{0})-\left\langle tu,u^{\ast}\right\rangle ]=f(x_{0}),
\]
and so $f(x_{0}+tu)=f(x_{0})+\left\langle tu,u^{\ast}\right\rangle $ for every
$t\in\mathbb{R}$. Hence (\ref{r-a4}) holds.

Taking $u\neq0$ and $u^{\ast}=0$, from (\ref{r-a4}) we have that
\begin{equation}
f_{\infty}(\pm u)=0\Longleftrightarrow\left[  f(x_{0}+tu)=f(x_{0})\ \ \forall
t\in\mathbb{R}\right]  \Leftrightarrow f|_{x_{0}+\mathbb{R}u}\text{ is
constant.} \label{r-a4a}%
\end{equation}

Moreover, it is worth observing that $f_{\infty}\geq0$ if $f$ is bounded from
below; indeed, if $f_{\infty}(u)<0$, from (\ref{r-a0}) we have that
$f(x+tu)\leq f(x)+tf_{\infty}(u)$, and so $\lim_{t\rightarrow\infty
}f(x+tu)=-\infty$, for every $x\in\operatorname*{dom}f$.

\smallskip In the sequel, for $\varphi,\psi:E\rightarrow\overline{\mathbb{R}}$
and $\rho\in\{\leq,<,=\}$ we set $[\varphi\;\rho\;\psi]:=\{x\in E\mid
\varphi(x)\;\rho\;\psi(x)\}$. For example $[\varphi\leq0]:=\{x\in X\mid
\varphi(x)\leq0\}$.

\medskip As in \cite[Def.\ 3]{Aza:19}, we say that $f$ is \emph{directionally
coercive} if $\lim_{t\rightarrow\infty}f(x+tu)=\infty$ for all $x\in X$ and
$u\in X\setminus\{0\}$, and $f$ is \emph{essentially directionally coercive}
if $f-x^{\ast}$ is directionally coercive for some $x^{\ast}\in X^{\ast}$.

From the equivalence of assertions (a'), (b') and (c') of Lemma
\ref{lem-Fact a} we get the next result.

\begin{corollary}
\label{cor-Fact b}Let $f\in\Gamma(X)$; then \emph{(a)} $f$ is directionally
coercive if and only if $[f_{\infty}\leq0]=\{0\}$, and \emph{(b)} $f$ is
essentially directionally coercive if and only if there exists $x^{\ast}\in
X^{\ast}$ such that $[f_{\infty}\leq x^{\ast}]=\{0\}$.
\end{corollary}

The previous result motivates a deeper study of proper lsc sublinear
functions; several properties of such functions are mentioned in
\cite[Th.\ 2.4.14]{Zal:02}.

Recall that the orthogonal spaces of the nonempty subsets $A\subset X$ and
$B\subset X^{\ast}$ are defined by
\[
A^{\perp}:=\{x^{\ast}\in X^{\ast}\mid\left\langle x,x^{\ast}\right\rangle
=0\ \forall x\in A\}\quad\text{and}\quad B^{\perp}:=\{x\in X\mid\left\langle
x,x^{\ast}\right\rangle =0\ \forall x^{\ast}\in B\},
\]
respectively; clearly, $A^{\perp}$ is a $w^{\ast}$-closed linear subspace of
$X^{\ast}$, $B^{\perp}$ is closed linear subspace of $X$, $A^{\perp}=\left(
\overline{\operatorname*{span}}A\right)  ^{\perp}$, $B^{\perp}=\left(
\overline{\operatorname*{span}}B\right)  ^{\perp}$, $(A^{\perp})^{\perp
}=\overline{\operatorname*{span}}A$, $(B^{\perp})^{\perp}=\overline
{\operatorname*{span}}B$. Also recall that the \emph{quasi-interior} and the
\emph{quasi-relative interior} of the convex subset $A$ of $X$ are the sets
\[
\operatorname*{qi}A:=\{a\in A\mid\overline{\mathbb{R}_{+}(A-a)}=X\},\quad
\operatorname*{qri}A:=\{a\in A\mid\overline{\mathbb{R}_{+}(A-a)}\text{ is a
linear space}\},
\]
respectively. Having in view that for $A\subset X$ a nonempty convex set one
has
\[
\mathbb{R}_{+}(A-a)\subset\operatorname*{span}(A-a)=\operatorname*{span}%
(A-A)=\mathbb{R}_{+}(A-A)\quad\forall a\in A, \label{r-a7}%
\]
one obtains (see e.g.\ \cite{Zal:15}) that
\begin{gather}
\operatorname*{qri}A=\{a\in A\mid\overline{\mathbb{R}_{+}(A-a)}=\overline
{\operatorname*{span}}(A-A)\}=\big\{a\in A\mid\overline{\mathbb{R}_{+}%
(A-a)}=\overline{\mathbb{R}_{+}(A-A)}\big\},\label{r-a11}\\
\operatorname*{qri}A=A\cap\operatorname*{qri}\overline{A},\quad
\operatorname*{qi}A=\left\{
\begin{array}
[c]{ll}%
\operatorname*{qri}A & \text{if }\overline{\mathbb{R}_{+}(A-A)}=X,\\
\emptyset & \text{otherwise.}%
\end{array}
\right.  \label{r-a12}%
\end{gather}

\section{Some results related to sublinear functions}

Throughout this section $g\in\Gamma(X)$ is assumed to be sublinear.

\begin{lemma}
\label{lem-Fact 1}Let us set $K:=[g\leq0]$ and $L:=K\cap(-K)$. Then $K$ is a
closed convex cone and $L$ is a closed linear subspace of $X$. Moreover,
\begin{gather}
L=\{x\in X\mid g(x)=g(-x)=0\}=[\partial g(0)]^{\perp},\label{r-a1a}\\
g(x+u)=g(x)\quad\forall x\in X,\ \forall u\in L. \label{r-a1}%
\end{gather}

\end{lemma}

Proof. Because $g$ is a (proper) lsc sublinear function, $[g\leq0]$ is clearly
a closed convex cone. The set $L$ is a closed convex cone as the intersection
of (two) closed convex cones. Since $L=-L$, $L$ is also a linear subspace of
$X$.

Take $x\in L$; because $0=g\left(  x+(-x)\right)  \leq g(x)+g(-x)\leq0+0=0$,
we get $g(x)=0=g(-x)$, and so $L\subset\{x\in X\mid g(x)=g(-x)=0\}$. The
reverse inclusion being obvious, the first equality in (\ref{r-a1a}) holds.

Set $B:=\partial g(0)$. Taking into account the formula for $g$ from
(\ref{r-a5b}), for $x\in X$ one has
\[
x\in L\Leftrightarrow g(\pm x)\leq0\Leftrightarrow\lbrack\pm\left\langle
x,x^{\ast}\right\rangle \leq0\ \forall x^{\ast}\in B]\Leftrightarrow\left[
\left\langle x,x^{\ast}\right\rangle =0\ \forall x^{\ast}\in B\right]
\Leftrightarrow x\in B^{\perp},
\]
and so the second equality in (\ref{r-a1a}) holds, too.

Take now $x\in X$ and $u\in L$. Using the sublinearity of $g$ one has
\[
g(x+u)\leq g(x)+g(u)=g(x)=g\left(  (x+u)+(-u)\right)  \leq
g(x+u)+g(-u)=g(x+u),
\]
and so $g(x+u)=g(x)$. \hfill$\square$

\begin{proposition}
\label{p-Fact 2}For $x^{\ast}\in X^{\ast}$ set $L_{x^{\ast}}:=\{x\in X\mid
g(\pm x)=\left\langle \pm x,x^{\ast}\right\rangle $. The following assertions hold:

\emph{(a)} If $x^{\ast}\in X^{\ast}$, then $L_{x^{\ast}}$ is a closed linear
subspace of $X$, and
\begin{gather}
L_{x^{\ast}}=\{x\in X\mid g(\pm x)\leq\left\langle \pm x,x^{\ast}\right\rangle
\}=\left[  \partial g(0)-x^{\ast}\right]  ^{\perp},\label{r-a1c}\\
g(x+u)=g(x)+\left\langle u,x^{\ast}\right\rangle \quad\forall x\in X,\ \forall
u\in L_{x^{\ast}}. \label{r-a1b}%
\end{gather}

\emph{(b)} If $u^{\ast}\in\partial g(0)$, then $L_{u^{\ast}}=\left[  \partial
g(0)-\partial g(0)\right]  ^{\perp}$. Consequently, $L_{x^{\ast}}\subset
L_{u^{\ast}}$ for all $x^{\ast}\in X^{\ast}$ and $u^{\ast}\in\partial g(0)$;
in particular $L_{u^{\ast}}=L_{v^{\ast}}$ for all $u^{\ast},v^{\ast}%
\in\partial g(0)$.
\end{proposition}

Proof. (a) Clearly, $h:=g-x^{\ast}$ is a proper lsc sublinear function. Using
Lemma \ref{lem-Fact 1} for $g$ replaced by $h$ we obtain that $L_{x^{\ast}}$
is a closed linear subspace of $X$ and the formulas for $L_{x^{\ast}}$ hold by
the definition of $L$ and because $\partial h(0)=\partial g(0)-x^{\ast}$.
Moreover, $g(x+u)-\left\langle x+u,x^{\ast}\right\rangle =g(x)-\left\langle
x,x^{\ast}\right\rangle $ for all $x\in X$ and $u\in L_{x^{\ast}}$, and so
(\ref{r-a1b}) holds, too.

(b) Take now $u^{\ast}\in\partial g(0)=:B$. Then $B-u^{\ast}\subset B-B$,
whence $Y:=\operatorname*{span}(B-u^{\ast})\subset\operatorname*{span}%
(B-B)=:Z$. Since $B-B=(B-u^{\ast})-(B-u^{\ast})\subset Y$, we get $Z\subset
Y$, and so $Y=Z$. Using (a) one has $L_{u^{\ast}}=B^{\perp}=Y^{\perp}%
=Z^{\perp}=(B-B)^{\perp}=\left[  \partial g(0)-\partial g(0)\right]  ^{\perp}$.

Let $x^{\ast}\in X^{\ast}$. Because $B-x^{\ast}\subset\operatorname*{span}%
(B-x^{\ast})$, as above one has $B-B\subset\operatorname*{span}(B-x^{\ast})$,
and so $L_{x^{\ast}}=(B-x^{\ast})^{\perp}=\left[  \operatorname*{span}%
(B-x^{\ast})\right]  ^{\perp}\subset(B-B)^{\perp}=L_{u^{\ast}}$.
\hfill$\square$

\medskip As seen in Proposition \ref{p-Fact 2}~(b), the set $\{L_{u^{\ast}%
}\mid u^{\ast}\in\partial g(0)\}$ is a singleton; its element will be denoted
by $L_{g}$ in the sequel.

\begin{proposition}
\label{p-Fact 2a}Let $x^{\ast}\in X^{\ast}$. The following assertions are
equivalent: \emph{(a)}$~x^{\ast}\in\operatorname*{qri}\partial g(0)$;
\emph{(b)}$~[g\leq x^{\ast}]$ is a linear space; \emph{(b')}~$L_{x^{\ast}%
}=[g\leq x^{\ast}]$; \emph{(c)}~$x^{\ast}\in\partial g(0)$ and $[g=x^{\ast}]$
is a linear space; \emph{(c')}~$x^{\ast}\in\partial g(0)$ and $L_{x^{\ast}%
}=[g=x^{\ast}]$.
\end{proposition}

Proof. Because $x^{\ast}\in\operatorname*{qri}\partial g(0)$ if and only if
$0\in\operatorname*{qri}\left[  \partial g(0)-x^{\ast}\right]  $ and
$\partial(g-x^{\ast})(0)=\partial g(0)-x^{\ast}$, we may (and do) assume that
$x^{\ast}=0$. Let us set $B:=\partial g(0)$ and $K:=[g\leq0]$; $K$ is a
(closed) convex cone and $l(K):=K\cap(-K)$ is a linear space.

Because $L_{x^{\ast}}$ is a linear space, the equivalences (b')
$\Leftrightarrow$ (b) and (c') $\Leftrightarrow$ (c) follow immediately from
(\ref{r-a1c}).

(c) $\Rightarrow$ (b) Because $0\in\partial g(0)$, one has $g\geq0$, and so
$[g\leq0]=[g=0]$. Hence (b) holds.

(b) $\Rightarrow$ (c) Because $K$ $(=[g\leq0]$) is a linear space, taking
$x\in K$ $(=-K)$ we get $g(\pm x)=0$ by Lemma \ref{lem-Fact 1}. It follows
that $g\geq0$ $(\Leftrightarrow0\in B)$ and $K=[g=0]$. Hence $[g=0]$ is a
linear space.

(b) $\Rightarrow$ (a) We have to show that $\overline{\mathbb{R}_{+}%
(B-B)}\subset\overline{\mathbb{R}_{+}B}$, the converse inclusion being
obvious. For this assume that $\overline{x}^{\ast}\in X^{\ast}\setminus
\overline{\mathbb{R}_{+}B}$. Then, by a separation theorem, there exist
$\overline{x}\in X$ and $\alpha\in\mathbb{R}$ such that $\left\langle
\overline{x},\overline{x}^{\ast}\right\rangle >\alpha\geq\left\langle
\overline{x},tu^{\ast}\right\rangle $ for all $t\in\mathbb{R}_{+}$ and
$u^{\ast}\in B$, whence $\alpha\geq0\geq\left\langle \overline{x},u^{\ast
}\right\rangle $ for $u^{\ast}\in B$, that is $\alpha\geq0\geq g(\overline
{x})$. Hence $0\neq\overline{x}\in K$ $(=-K)$, and so $g(\pm\overline{x})=0$.
It follows that $\left\langle \pm x,u^{\ast}\right\rangle \leq g(\pm x)=0$,
whence $\left\langle x,u^{\ast}\right\rangle =0$, for all $u^{\ast}\in B$.
Hence $\left\langle \overline{x},t(u^{\ast}-v^{\ast})\right\rangle
=0<\left\langle \overline{x},\overline{x}^{\ast}\right\rangle $ for all
$t\in\mathbb{R}_{+}$ and $u^{\ast},v^{\ast}\in B$, proving that $\overline
{x}\notin\overline{\mathbb{R}_{+}(B-B)}$. Therefore, $(x^{\ast}=)$
$0\in\operatorname*{qri}B$.

(a) $\Rightarrow$ (b) Because $0\in\operatorname*{qri}B$, $0\in B$ and
$\operatorname*{cl}(\mathbb{R}_{+}B)=\overline{\mathbb{R}_{+}(B-B)}$. Take
$x\in K$; then $\left\langle x,u^{\ast}\right\rangle \leq g(x)\leq0$ for
$u^{\ast}\in B$, and so $\left\langle -x,u^{\ast}\right\rangle \geq0$ for all
$u^{\ast}\in B$, whence $-x\in\left[  \operatorname*{cl}(\mathbb{R}%
_{+}B)\right]  ^{+}=\big(  \overline{\mathbb{R}_{+}(B-B)}\big)
^{+}$. It follows that $\left\langle x,v^{\ast}\right\rangle =\left\langle
-x,0-v^{\ast}\right\rangle \geq0$, that is $\left\langle -x,v^{\ast
}\right\rangle \leq0$, for all $v^{\ast}\in B$, whence $g(-x)\leq0$. Hence
$x\in-K$, and so $K$ is a linear space. \hfill$\square$

\begin{corollary}
\label{cor-Fact 2b}Let $x^{\ast}\in X^{\ast}$. Then $x^{\ast}\in
\operatorname*{qi}\partial g(0)$ if and only if $[g\leq x^{\ast}]=\{0\}$.
\end{corollary}

Proof. Set $B:=\partial g(0).$ Assume that $x^{\ast}\in\operatorname*{qi}B$.
From (\ref{r-a11}) and (\ref{r-a12}) we have that $x^{\ast}\in
\operatorname*{qri}B$ and $\overline{\mathbb{R}_{+}(B-B)}=X^{\ast}$. Using the
equivalence (a)~$\Leftrightarrow$ (b') of Proposition \ref{p-Fact 2a} and
Proposition \ref{p-Fact 2}~(b) we obtain that $[g\leq x^{\ast}]=L_{x^{\ast}%
}=(X^{\ast})^{\perp}=\{0\}$.

Conversely, assume that $[g\leq x^{\ast}]=\{0\}$. Using the implication
(b)~$\Rightarrow$ (a)~$\wedge$ (b') of Proposition \ref{p-Fact 2a}, we get
$x^{\ast}\in\operatorname*{qri}B$ $[\subset B]$ and $(L_{g}=)$ $L_{x^{\ast}%
}=\{0\}$. Using now Proposition \ref{p-Fact 2}~(b) we obtain that $X^{\ast
}=\{0\}^{\perp}=L_{x^{\ast}}^{\perp}=\big(\left[  B-B\right]  ^{\perp
}\big)^{\perp}=\overline{\mathbb{R}_{+}(B-B)}$. Using again (\ref{r-a12}) we
get $x^{\ast}\in\operatorname*{qi}B$. \hfill$\square$

\begin{proposition}
\label{p-Fact 2c}Assume that $X$ is a separable normed vector space. Then
$w^{\ast}$-$\operatorname*{qri}\partial g(0)\neq\emptyset$, and so there
exists $x^{\ast}\in X^{\ast}$ such that the set $[g\leq x^{\ast}]$ is a linear space.
\end{proposition}

Proof. In order to get the conclusion we apply \cite[Th.\ 2.19~(b)]{BorLew:92}
which states that for any weakly$^{\ast}$ cs-closed subset $C\subset(X^{\ast
},w^{\ast})$, $X$ being a separable nvs, one has $w^{\ast}$%
-$\operatorname*{qri}C\neq\emptyset$. So, consider $(\alpha_{n})_{n\geq
1}\subset\mathbb{R}_{+}$ and $(x_{n}^{\ast})_{n\geq1}\subset C:=\partial g(0)$
such that $w^{\ast}$-$\lim\sum_{k=1}^{n}\alpha_{n}x_{n}^{\ast}=x^{\ast}\in
X^{\ast}$. We need to prove that $x^{\ast}\in C$. For this, observe first that
there exists $n_{0}\geq1$ such that $\alpha_{n_{0}}>0$. Then, for $n\geq
n_{0}$ we have that $\beta_{n}:=\sum_{k=1}^{n}\alpha_{n}>0$ and $u_{n}^{\ast
}:=\beta_{n}^{-1}\sum_{k=1}^{n}\alpha_{n}x_{n}^{\ast}\in C$. Since $\beta
_{n}\rightarrow1$, we obtain that $C\ni w^{\ast}$-$\lim u_{n}^{\ast}=x^{\ast}%
$. The proof is complete. \hfill$\square$

\begin{remark}
\label{rem2}Notice that the separability of the nvs $X$ in Proposition
\ref{p-Fact 2c} is essential. For example, the space of square summable
real-valued functions $X:=\ell_{2}(\Gamma)$, endowed with the norm $\left\Vert
\cdot\right\Vert $ defined by $\left\Vert x\right\Vert :=\big(  \sum
_{\gamma\in\Gamma}\left\vert x(\gamma)\right\vert ^{2}\big)  ^{1/2}$, is a
Hilbert space, while $X_{+}:=\{x\in X\mid x(\gamma)\geq0$ $\forall\gamma
\in\Gamma\}$ is a closed convex cone such that $X_{+}-X_{+}=X$. If $\Gamma$ is
at most countable, then $\operatorname*{qri}X_{+}=\operatorname*{qi}%
X_{+}=\{x\in X\mid x(\gamma)>0$ $\forall\gamma\in\Gamma\}$. If $\Gamma$ is
uncountable then, as in \cite[Ex.\ 3.11~(iii)]{BorLew:92},
$\operatorname*{qri}X_{+}=\emptyset$.
\end{remark}

\medskip Considering the quotient space $\widehat{X}:=X/L_{g}:=\{\widehat
{x}\mid x\in X\}$ of $X$ with respect to $L_{g}$ endowed with the quotient
topology, $\widehat{X}$ becomes a separated locally convex space such that the
natural projection $\pi:X\rightarrow\widehat{X}$, defined by $\pi
(x):=\widehat{x}$, is a continuous open linear operator; moreover
$A\subset\widehat{X}$ is closed if and only if $\pi^{-1}(A)$ is closed.

Fixing $x^{\ast}\in\partial g(0)$ one has $L_{x^{\ast}}=L_{g}$; using
(\ref{r-a1b}), we obtain that
\begin{equation}
\widehat{g}_{x^{\ast}}:\widehat{X}\rightarrow\overline{\mathbb{R}}%
,\quad\widehat{g}_{x^{\ast}}(\widehat{x}):=g(x)-\left\langle x,x^{\ast
}\right\rangle \quad(x\in X) \label{r-gxs}%
\end{equation}
is well defined.

\begin{proposition}
\label{p-Fact 3}Assume that $x^{\ast}\in\partial g(0)$. Then $\widehat
{g}_{x^{\ast}}$ defined by (\ref{r-gxs}) is a proper lsc sublinear function
such that $\widehat{g}_{x^{\ast}}\geq0\ $and $L_{\widehat{g}_{x^{\ast}}%
}=\{\widehat{0}\}$. Moreover, $x^{\ast}\in\operatorname*{qri}\partial g(0)$ if
and only if $0\in\operatorname*{qi}\partial\widehat{g}_{x^{\ast}}(\widehat
{0})$.
\end{proposition}

Proof. The fact that $\widehat{g}_{x^{\ast}}$ is proper, sublinear and takes
nonnegative values follows immediately from its definition. For $\alpha
\in\mathbb{R}$ one has
\begin{equation}
\lbrack\widehat{g}_{x^{\ast}}\leq\alpha]=\{\widehat{x}\in\widehat{X}%
\mid\widehat{g}_{x^{\ast}}(\widehat{x})\leq\alpha\}=\pi(\{x\in X\mid
g(x)-\left\langle x,x^{\ast}\right\rangle \leq\alpha\})=\pi([g-x^{\ast}%
\leq\alpha]); \label{r-a8}%
\end{equation}
using (\ref{r-a1}) we have that $\pi^{-1}([\widehat{g}_{x^{\ast}}\leq
\alpha])=[g-x^{\ast}\leq\alpha]$. Since $g$ is lsc, $g-x^{\ast}$ is so; it
follows that $[g-x^{\ast}\leq\alpha]$ is closed, and so $[\widehat{g}%
_{x^{\ast}}\leq\alpha]$ is closed in $\widehat{X}$ for every $\alpha
\in\mathbb{R}$, whence $\widehat{g}_{x^{\ast}}$ is lsc.

Because $\widehat{g}_{x^{\ast}}\geq0$ one has $0\in\partial\widehat
{g}_{x^{\ast}}(0)$, and so $L_{\widehat{g}_{x^{\ast}}}=\{\widehat{x}%
\mid\widehat{g}_{x^{\ast}}(\widehat{x})=\widehat{g}_{x^{\ast}}(-\widehat
{x})=0\}$. Take $x\in X$ with $\widehat{x}\in L_{\widehat{g}_{x^{\ast}}}$;
from the definition of $\widehat{g}_{x^{\ast}}$ we have that $g(\pm
x)-\left\langle \pm x,x^{\ast}\right\rangle =0$, and so $x\in L_{x^{\ast}%
}=L_{g}$. It follows that $\widehat{x}=0$, and so $L_{\widehat{g}_{x^{\ast}}%
}=\{\widehat{0}\}$. Taking $\alpha:=0$ in (\ref{r-a8}) and in the equality on
the line below it we obtain that $[\widehat{g}_{x^{\ast}}\leq0]=\pi
(K_{x^{\ast}})$ and $\pi^{-1}([\widehat{g}_{x^{\ast}}\leq0])=K_{x^{\ast}}$;
hence $K_{x^{\ast}}$ is a linear space if and only if $[\widehat{g}_{x^{\ast}%
}\leq0]$ is a linear space. Using Proposition \ref{p-Fact 2a} we obtain that
$x^{\ast}\in\operatorname*{qri}\partial g(0)$ if and only if $0\in
\operatorname*{qri}\widehat{g}_{x^{\ast}}(0)$. Because $L_{\widehat
{g}_{x^{\ast}}}=\{\widehat{0}\}$ we have that $\operatorname*{cl}%
\big[\mathbb{R}_{+}\big(\partial\widehat{g}_{x^{\ast}}(\widehat{0}%
)-\partial\widehat{g}_{x^{\ast}}(\widehat{0})\big)\big]=(\widehat{X})^{\ast}$,
and so $\operatorname*{qri}\widehat{g}_{x^{\ast}}(0)=\operatorname*{qi}%
\widehat{g}_{x^{\ast}}(0)$ by (\ref{r-a12}). \hfill$\square$

\medskip

In this context it is natural to know sufficient conditions for having
$g(x)>0$ for $x\in X\setminus\{0\}$. Some sufficient conditions are provided
in the next result. Recall that the \emph{core} of the subset $A$ of the real
linear space $E$ is $\operatorname*{core}A:=\{x\in E\mid\forall u\in E,$
$\exists\delta>0,$ $\forall t\in\lbrack0,\delta]:x+tu\in A\}.$

\begin{proposition}
\label{p-Fact 3b}Let $x^{\ast}\in X^{\ast}$. Consider the following assertions:

\emph{(i)}$~[g\leq0]=\{0\}$;

\emph{(ii)} $x^{\ast}\in\operatorname*{qi}\partial g(0)$;

\emph{(iii)}~$x^{\ast}\in\operatorname*{core}\partial g(0)$;

\emph{(iv)} $x^{\ast}\in\operatorname*{int}_{\tau}\partial g(0)$, where $\tau$
is a linear topology on $X^{\ast}$;

\emph{(v)} the topology of $X$ is defined by the norm $\left\Vert
\cdot\right\Vert $ and $x^{\ast}\in\operatorname*{int}_{\left\Vert
\cdot\right\Vert _{\ast}}\left(  \partial g(0)\right)  $, where $\left\Vert
\cdot\right\Vert _{\ast}$ is the dual norm on $X^{\ast}$;

\emph{(vi)} the topology of $X$ is defined by the norm $\left\Vert
\cdot\right\Vert $ and there exists $\alpha>0$ such that $g(x)-\left\langle
x,x^{\ast}\right\rangle \geq\alpha\left\Vert x\right\Vert $ for all $x\in X$.

Then \emph{(vi)}~$\Leftrightarrow$ \emph{(v)}~$\Rightarrow$ \emph{(iv)}%
~$\Rightarrow$ \emph{(iii)}~$\Rightarrow$ \emph{(ii)}~$\Leftrightarrow$
\emph{(i)}; moreover, if $\dim X<\infty$ then \emph{(i)}~$\Rightarrow$
\emph{(vi)}.
\end{proposition}

Proof. Because $\partial\left(  g-x^{\ast}\right)  (0)=\partial g(0)-x^{\ast}%
$, we may (and do) assume that $x^{\ast}=0$. We set $B:=\partial g(0)$.

(vi)~$\Leftrightarrow$ (v) This assertion follows immediately from the
equivalence of assertions (e) and (f) of \cite[Exer.\ 2.41]{Zal:02}.

(v)~$\Rightarrow$ (iv) This assertion is true because the topology generated
by any norm on a linear space is a linear topology.

(iv)~$\Rightarrow$ (iii) It is well known that $\operatorname*{core}%
A=\operatorname*{int}_{\sigma}A$ when $A$ is a convex subset of topological
vector space $(Y,\sigma)$ with $\operatorname*{int}_{\sigma}A\neq\emptyset$.
The set $B\subset X^{\ast}$ being convex, the implication is true.

(iii)~$\Rightarrow$ (ii) Because $0\in\operatorname*{core}B$, we have that
$0\in B$ and $\mathbb{R}_{+}B=X^{\ast}$, and so $\operatorname*{cl}%
(\mathbb{R}_{+}B)=X^{\ast}$. Therefore, $0\in\operatorname*{qi}B$.

(ii)~$\Leftrightarrow$ (i) This equivalence is provided by Corollary
\ref{cor-Fact 2b}.

(i)~$\Rightarrow$ (vi) (if $\dim X<\infty$). Assume that $\dim X<\infty$. It
is well known that all the norms on a finite dimensional linear space are
equivalent, and a separated linear topology on such a space is normable. So,
let let $\left\Vert \cdot\right\Vert $ be a norm on $X$. Because
$S_{X}:=\{x\in X\mid\left\Vert x\right\Vert =1\}$ is compact and $g$ is lsc,
there exists $\overline{x}\in S_{X}$ such that $g(x)\geq g(\overline
{x})=:\alpha$ $(>0)$. Taking $x\in X\setminus\{0\}$, $x^{\prime}:=\left\Vert
x\right\Vert ^{-1}x\in S_{X}$, and so $g(x)=\left\Vert x\right\Vert \cdot
g(x^{\prime})\geq\alpha\left\Vert x\right\Vert $. Hence (vi) holds.
\hfill$\square$

\medskip Note that the reverse implication of (i)~$\Rightarrow$ (vi) from
Proposition \ref{p-Fact 3b} is not true even when $X$ is an infinite
dimensional separable Hilbert space. Indeed, take $X:=\ell_{2}$ endowed with
its usual norm $\left\Vert \cdot\right\Vert _{2}$ and $g:X\rightarrow
\overline{\mathbb{R}}$ defined by $g(x):=\big(\sum_{n\geq1}\left\vert
x_{n}\right\vert ^{q}\big)^{1/q}$ $(=\left\Vert x\right\Vert _{q})$ for
$x:=(x_{n})_{n\geq1}\in X$, where $q\in{}]2,\infty\lbrack$. Because $\ell
_{p}\subset\ell_{p^{\prime}}$ for $1\leq p<p^{\prime}\leq\infty$ with
$\left\Vert x\right\Vert _{p^{\prime}}\leq\left\Vert x\right\Vert _{p}$ for
$x\in\ell_{p}$, $g(x)\leq\left\Vert x\right\Vert _{2}$ for all $x\in X$, and
so $g$ is a finitely valued continuous sublinear function verifying (i).
Assuming that (vi) holds, there exists $\alpha>0$ such that $\left\Vert
x\right\Vert _{q}\geq\alpha\left\Vert x\right\Vert _{2}$ for all $x\in\ell
_{2}$. Consider the sequence $x=(n^{-1/2})_{n\geq1}\subset\mathbb{R}$; then
$\xi_{n}:=(1,...,n^{-1/2},0,0,...)\in\ell_{2}$, and so
\[
\left(  \sum\nolimits_{k=1}^{n}\frac{1}{k}\right)  ^{1/2}=\left\Vert \xi
_{n}\right\Vert _{2}\leq\alpha^{-1}\left\Vert \xi_{n}\right\Vert _{q}%
=\alpha^{-1}\left(  \sum\nolimits_{k=1}^{n}\frac{1}{k^{q/2}}\right)
^{1/q}\quad\forall n\geq1,
\]
whence the contradiction $\infty=\lim_{n\rightarrow\infty}\left(  \sum
_{k=1}^{n}\frac{1}{k}\right)  ^{1/2}\leq\lim_{n\rightarrow\infty}\alpha
^{-1}\big(\sum_{k=1}^{n}\frac{1}{k^{q/2}}\big)^{1/q}<\infty$.

\section{Applications to the shape of convex functions}

The following results are motivated by the notion and results from
\cite{Aza:19}.

In the sequel, for $f\in\Gamma(X)$ we set $L_{f}:=L_{f_{\infty}}$. As seen in
(\ref{r-a5}), $\partial f_{\infty}(0)=\overline{\operatorname*{dom}f^{\ast}}$,
and so, by Proposition \ref{p-Fact 2}, we have that
\begin{equation}
L_{f}=\{u\in X\mid f_{\infty}(\pm u)=\left\langle \pm u,x^{\ast}\right\rangle
\}\text{ for some (any) }x^{\ast}\in\overline{\operatorname*{dom}f^{\ast}}.
\label{r-a15}%
\end{equation}
Also, observe that $\operatorname{Im}\partial f\subset\operatorname*{dom}%
f^{\ast}\subset\overline{\operatorname*{dom}f^{\ast}}$ and
\[
\overline{\operatorname*{span}}\left(  \operatorname{Im}\partial
f-\operatorname{Im}\partial f\right)  \subset\overline{\operatorname*{span}%
}\left(  \operatorname*{dom}f^{\ast}-\operatorname*{dom}f^{\ast}\right)
=\overline{\operatorname*{span}}\left(  \overline{\operatorname*{dom}f^{\ast}%
}-\overline{\operatorname*{dom}f^{\ast}}\right)  =(L_{f})^{\perp},
\]
the first inclusion becoming equality if $X$ is a Banach space because in this
case $\operatorname*{dom}f^{\ast}\subset\operatorname*{cl}_{\left\Vert
\cdot\right\Vert }\operatorname{Im}\partial f\subset\overline
{\operatorname*{dom}f^{\ast}}$ by Br{\o }ndsted--Rockafellar theorem (see
e.g.\ \cite[Th.\ 3.1.2]{Zal:02}).

\begin{corollary}
\label{cor-Fact 4a}\emph{(a)} The function $f$ is directionally coercive if
and only if $0\in\operatorname*{qi}\overline{\operatorname*{dom}f^{\ast}}$.

\emph{(b)} The function $f$ is essentially directionally coercive if and only
if $\operatorname*{qi}\overline{\operatorname*{dom}f^{\ast}}\neq\emptyset$.
\end{corollary}

Proof. Having $x^{\ast}\in X^{\ast}$, one has $(f-x^{\ast})_{\infty}%
=f_{\infty}-x^{\ast}$ and $\operatorname*{dom}(f-x^{\ast})^{\ast
}=\operatorname*{dom}f^{\ast}-x^{\ast}$, and so $\operatorname*{qi}%
\overline{\operatorname*{dom}(f-x^{\ast})}=\operatorname*{qi}\overline
{\operatorname*{dom}f^{\ast}}-x^{\ast}$. Hence (a)~$\Rightarrow$ (b) by
Corollary \ref{cor-Fact b}.

(a) By Corollary \ref{cor-Fact b} one has that $f$ is directionally coercive
if and only $[f_{\infty}\leq0]=\{0\}$, and the latter is equivalent to
$0\in\operatorname*{qi}\overline{\operatorname*{dom}f^{\ast}}$ by (\ref{r-a5})
and the equivalence (i)~$\Leftrightarrow$ (ii) of Proposition \ref{p-Fact 3b}.
\hfill$\square$

\medskip The representation of the (continuous) convex function $f$ from
Theorems 4--5 of \cite{Aza:19} motivates the next result.

\begin{proposition}
\label{p-Fact 4}Assume that $f=h\circ A+x^{\ast}$, where $x^{\ast}\in X^{\ast
}$, $h\in\Gamma(Y)$ with $Y$ a separated locally convex space, and
$A:X\rightarrow Y$ is a continuous linear operator. Then the following
assertions hold:

\emph{(a)}$~f_{\infty}=h_{\infty}\circ A+x^{\ast}$ and $\ker A\subset\lbrack
f_{\infty}=x^{\ast}]$;

\emph{(b)}~if $h_{\infty}\geq0$, then $x^{\ast}\in\partial f_{\infty
}(0)=\overline{\operatorname*{dom}f^{\ast}}$, the converse implication being
true if, moreover, $\operatorname{Im}A=Y$;

\emph{(c)}~if $[h_{\infty}\leq0]=\{0\}$ then $\ker A=[f_{\infty}\leq x^{\ast
}]$ and $x^{\ast}\in\operatorname*{qri}\overline{\operatorname*{dom}f^{\ast}}%
$; conversely, if $\operatorname{Im}A=Y$ and $\ker A=[f_{\infty}\leq x^{\ast
}]$, then $[h_{\infty}\leq0]=\{0\}$.

\emph{(d)} if $h$ is bounded from below, then $h_{\infty}\geq0$ and $x^{\ast
}\in\operatorname*{dom}f^{\ast}$; if $\operatorname{Im}A=Y$ then $\inf
h=-f^{\ast}(x^{\ast})$, and so $h$ is bounded from below if and only if
$x^{\ast}\in\operatorname*{dom}f^{\ast}$.

\emph{(e)} Assume that $\operatorname{Im}A=Y$. Then $h$ attains its infimum on
$Y$ if and only if $x^{\ast}\in\operatorname{Im}\partial f$.
\end{proposition}

\smallskip Proof. (a) Let $x_{0}\in\operatorname*{dom}f$; then $Ax_{0}%
\in\operatorname*{dom}h$ and
\begin{align*}
f_{\infty}(u)  &  =\lim_{t\rightarrow\infty}\frac{f(x_{0}+tu)-f(x_{0})}%
{t}=\lim_{t\rightarrow\infty}\frac{h(Ax_{0}+tAu)-h(Ax_{0})+t\left\langle
u,x^{\ast}\right\rangle }{t}\\
&  =h_{\infty}(Au)+\left\langle u,x^{\ast}\right\rangle \quad\forall u\in X,
\end{align*}
and so $f_{\infty}=h_{\infty}\circ A+x^{\ast}$. The desired inclusion follows
now immediately.

(b) Assume that $h_{\infty}\geq0$. Then $f_{\infty}-x^{\ast}=h_{\infty}\circ
A\geq0$, and so $x^{\ast}\in\partial f_{\infty}(0)=\overline
{\operatorname*{dom}f^{\ast}}$. Assume now that $x^{\ast}\in\overline
{\operatorname*{dom}f^{\ast}}$ and $\operatorname{Im}A=Y$. Clearly $x^{\ast
}\in\partial f_{\infty}(0)$, and so $f_{\infty}\geq x^{\ast}$. Taking $y\in
Y=\operatorname{Im}A$, there exists $u\in X$ with $Au=y$. Hence $h_{\infty
}(y)=h_{\infty}(Au)=f_{\infty}(u)-\left\langle u,x^{\ast}\right\rangle \geq0$.
Therefore $h_{\infty}\geq0$.

(c) Assume that $[h_{\infty}\leq0]=\{0\}$. Then $h_{\infty}\geq0$, and so
$x^{\ast}\in\overline{\operatorname*{dom}f^{\ast}}$ by (b); moreover, by (a),
$\ker A\subset\lbrack f_{\infty}=x^{\ast}]$. Take $u\in\lbrack f_{\infty
}=x^{\ast}]$; then $h_{\infty}(Au)=f_{\infty}(u)-\left\langle u,x^{\ast
}\right\rangle =0$, and so $Au=0$, that is $u\in\ker A$. Hence $\ker
A=[f_{\infty}=x^{\ast}]$; this shows that $[f_{\infty}=x^{\ast}]$ is a linear
space and so, using the implications (c)~$\Rightarrow$ (b)~$\Rightarrow$ (a)
of Proposition \ref{p-Fact 2a}, we obtain that $\ker A=[f_{\infty}\leq
x^{\ast}]$ and $x^{\ast}\in\operatorname*{qri}\overline{\operatorname*{dom}%
f^{\ast}}$.

Assume now that $\operatorname{Im}A=Y$ and $\ker A=[f_{\infty}\leq x^{\ast}]$.
Using the implication (b)~$\Rightarrow$ (c) of Proposition \ref{p-Fact 2a} we
obtain that $x^{\ast}\in\overline{\operatorname*{dom}f^{\ast}}$ and $\ker
A=[f^{\ast}=x^{\ast}]$. From (b) we have that $h_{\infty}\geq0$. Take
$y\in\lbrack h_{\infty}\leq0]$. Because $\operatorname{Im}A=Y$, there exists
$u\in X$ such that $y=Au$, and so $f_{\infty}(u)=h_{\infty}(Au)+\left\langle
u,x^{\ast}\right\rangle \leq\left\langle u,x^{\ast}\right\rangle $. Hence
$u\in\ker A$, and so $y=Au=0$.

(d) Assume that $h$ is bounded from below. Then $f-x^{\ast}\geq\inf
h\in\mathbb{R}$, and so, $f^{\ast}(x^{\ast})=\sup_{x\in X}[\left\langle
x,x^{\ast}\right\rangle -f(x)]<\infty$, whence $x^{\ast}\in\operatorname*{dom}%
f^{\ast}$.

Assume now that $\operatorname{Im}A=Y$. Then $\inf h=\inf h\circ
A=\inf(f-x^{\ast})=-f^{\ast}(x^{\ast})$. Hence, $h$ is bounded from below if
and only if $x^{\ast}\in\operatorname*{dom}f^{\ast}$.

(e) Assume $\operatorname{Im}A=Y$. Suppose that $h$ attains its infimum at
$\overline{y}\in Y$ and take $\overline{x}\in X$ such that $A\overline
{x}=\overline{y}$. Then $f(x)-\left\langle x,x^{\ast}\right\rangle =h(Ax)\geq
h(A\overline{x})=f(\overline{x})-\left\langle \overline{x},x^{\ast
}\right\rangle $, whence $x^{\ast}\in\partial f(\overline{x})\subset
\operatorname{Im}\partial f$. Conversely, assume that $x^{\ast}\in\partial
f(\overline{x})$ $(\subset\operatorname*{dom}f^{\ast})$. Then, as seen in (d),
$\inf h=-f^{\ast}(x^{\ast})=f(\overline{x})-\left\langle \overline{x},x^{\ast
}\right\rangle =h(A\overline{x})$. \hfill$\square$

\begin{lemma}
\label{lem-Fact 5}Let $x^{\ast}\in X^{\ast}$ and $L:=L_{x^{\ast}}:=\{u\in
X\mid f_{\infty}(\pm u)=\pm\left\langle u,x^{\ast}\right\rangle \}$. Then $L$
is a closed linear subspace of $X$, $\operatorname*{dom}%
f+L=\operatorname*{dom}f$, $(X\setminus\operatorname*{dom}f)+L=X\setminus
\operatorname*{dom}f$, and
\begin{equation}
f(x+u)=f(x)+\left\langle u,x^{\ast}\right\rangle \quad\forall x\in X,\ \forall
u\in L. \label{r-a2}%
\end{equation}

\end{lemma}

Proof. Applying Lemma \ref{lem-Fact 1} for $g:=f_{\infty}$, we have that $L$
is a closed linear subspace of $X$. Because $0\in L$ the inclusions
$\operatorname*{dom}f+L\supset\operatorname*{dom}f$ and $(X\setminus
\operatorname*{dom}f)+L\supset X\setminus\operatorname*{dom}f$ are obvious.
Take $x\in\operatorname*{dom}f$ and $u\in L$. Then $f(x+u)\leq f(x)+f_{\infty
}(u)=f(x)+\left\langle u,x^{\ast}\right\rangle <\infty$, and so
$\operatorname*{dom}f+L\subset\operatorname*{dom}f$; hence
$\operatorname*{dom}f+L=\operatorname*{dom}f$. Assuming that for some $x\in
X\setminus\operatorname*{dom}f$ and $u\in L$ one has $x^{\prime}%
:=x+u\in\operatorname*{dom}f$ we get the contradiction $x=x^{\prime}%
+(-u)\in\operatorname*{dom}f$. Hence $(X\setminus\operatorname*{dom}%
f)+L=X\setminus\operatorname*{dom}f$.

From the previous equality it is clear that $f(x+u)=f(x)+\left\langle
u,x^{\ast}\right\rangle $ $(=\infty)$ for $x\in X\setminus\operatorname*{dom}%
f$ and $u\in L$. Take now $x\in\operatorname*{dom}f$ and $u\in L$. Then
$x+u\in\operatorname*{dom}f$ and, as seen above, $f(x+u)\leq f(x)+\left\langle
u,x^{\ast}\right\rangle $. Hence
\[
f(x+u)\leq f(x)+\left\langle u,x^{\ast}\right\rangle \leq f(x+u)+\left\langle
-u,x^{\ast}\right\rangle +\left\langle u,x^{\ast}\right\rangle =f(x+u),
\]
and so $f(x+u)=f(x)+\left\langle u,x^{\ast}\right\rangle $. Therefore,
(\ref{r-a2}) holds. \hfill$\square$

\medskip

In the conditions and notation of Lemma \ref{lem-Fact 5} we have that
$f(x+u)-\left\langle x+u,x^{\ast}\right\rangle =f(x)-\left\langle x,x^{\ast
}\right\rangle $ for all $x\in X$ and $u\in L$, which shows that
\begin{equation}
h_{x^{\ast}}:X/L\rightarrow\overline{\mathbb{R}},\quad h_{x^{\ast}}%
(\widehat{x}):=f(x)-\left\langle x,x^{\ast}\right\rangle \quad(x\in X)
\label{r-a3}%
\end{equation}
is well defined and $f=h_{x^{\ast}}\circ\pi+x^{\ast}$, where $\pi:X\rightarrow
X/L$ is the (natural) projection defined by $\pi(x):=\widehat{x}$. The
convexity and properness of $h$ follow immediately from the corresponding
properties of $f$.

\begin{proposition}
\label{p-Fact 6}Let $x^{\ast}\in X^{\ast}$, $L:=L_{x^{\ast}}:=\{u\in X\mid
f_{\infty}(\pm u)=\pm\left\langle u,x^{\ast}\right\rangle \}$, and
$h:=h_{x^{\ast}}$ be defined in (\ref{r-a3}). Then the following assertions hold:

\emph{(a)} $h\in\Gamma(X/L)$, $h_{\infty}(\widehat{u})=f_{\infty
}(u)-\left\langle u,x^{\ast}\right\rangle $ for all $u\in X$, and
$\{\widehat{u}\in X/L\mid h_{\infty}(\widehat{u})=h_{\infty}(-\widehat
{u})=0\}=\{\widehat{0}\}$;

\emph{(b)} $h_{\infty}\geq0$ if and only if $x^{\ast}\in\overline
{\operatorname*{dom}f^{\ast}}$;

\emph{(c)} if $x^{\ast}\in\overline{\operatorname*{dom}f^{\ast}}$
(consequently $L=L_{f}$), then%
\[
\lbrack h_{\infty}\leq0]=\{\widehat{0}\}\Longleftrightarrow x^{\ast}%
\in\operatorname*{qri}\overline{\operatorname*{dom}f^{\ast}}%
\Longleftrightarrow L=[f_{\infty}\leq x^{\ast}]\Longleftrightarrow
L=[f_{\infty}=x^{\ast}];
\]

\emph{(d)} $\inf h=-f^{\ast}(x^{\ast})$, and so $h$ is bounded from below if
and only if $x^{\ast}\in\operatorname*{dom}f^{\ast}$;

\emph{(e)} $h$ attains its infimum on $X/L$ if and only if $x^{\ast}%
\in\operatorname{Im}\partial f$.
\end{proposition}

Proof. (a) As seen above, $h$ is well defined, proper and convex, and
$f=h\circ\pi+x^{\ast}$, where $\pi:X\rightarrow X/L$ with $\pi(x):=\widehat
{x}$. For $\alpha\in\mathbb{R}$ and $x\in X$ one has
\[
\widehat{x}\in\lbrack h\leq\alpha]\Leftrightarrow f(x)-\left\langle x,x^{\ast
}\right\rangle \leq\alpha\Leftrightarrow x\in\lbrack f-x^{\ast}\leq\alpha],
\]
and so $\pi^{-1}\left(  [h\leq\alpha]\right)  =[f-x^{\ast}\leq\alpha]$. Since
$f-x^{\ast}\in\Gamma(X)$, $[f-x^{\ast}\leq\alpha]$ is closed. Hence
$[h\leq\alpha]$ is closed. Because $\alpha\in\mathbb{R}$ is arbitrary, it
follows that $h$ is lsc. Therefore, $h\in\Gamma(X/L)$. The expression of
$h_{\infty}$ is obtained using Proposition \ref{p-Fact 4}~(a). Take $u\in X$;
from the expression of $h_{\infty}$ we have that
\[
h_{\infty}(\pm\widehat{u})=0\Longleftrightarrow h_{\infty}(\widehat{\pm
u})=0\Longleftrightarrow f_{\infty}(\pm u)-\left\langle \pm u,x^{\ast
}\right\rangle =0\Longleftrightarrow u\in L\Longleftrightarrow\widehat
{u}=\widehat{0}.
\]

Because $\pi$ is onto, the assertions (b), (d), (e) follow from assertions
(b), (d) and (e) of Proposition \ref{p-Fact 4}, respectively.

(c) Because $\pi$ is onto and $\ker\pi=L$, using Proposition \ref{p-Fact 4}%
~(c) we get the equivalence $[h_{\infty}\leq0]=\{\widehat{0}\}~\Leftrightarrow
$ $L=[f\leq x^{\ast}]$; the other equivalences follow from
(a)~$\Leftrightarrow$ (b')~$\Leftrightarrow$ (c') of Proposition
\ref{p-Fact 2a} because $x^{\ast}\in\overline{\operatorname*{dom}f^{\ast}%
}=\partial f_{\infty}(0)$. \hfill$\square$

\medskip Our main result is the following theorem; in its statement, for the
closed linear subspace $Y$ of $X$, $\pi:X\rightarrow X/Y$ is the natural
projection of $X$ onto $Y$, that is $\pi(x):=\widehat{x}$.

\begin{theorem}
\label{t-aza}Let $f\in\Gamma(X)$. The following assertions hold:

\emph{(i)} For every $x^{\ast}\in X^{\ast}$, there exist a closed linear
subspace $Y$ of $X$ and $h\in\Gamma(X/Y)$ such that $h$ is not constant on any
line $\widehat{x}+\mathbb{R}\widehat{u}$ with $\widehat{u}\neq\widehat{0}$
such that $f=h\circ\pi+x^{\ast}$. Moreover, for $x^{\ast}\in
\operatorname*{dom}f^{\ast}$, $h$ is bounded from below, while for $x^{\ast
}\in\operatorname{Im}\partial f$, $h$ attains its infimum on $X/Y$; in both
cases $Y=L_{f}$.

\emph{(ii)} There exist a closed linear subspace $Y$ of $X$, a directionally
coercive function $h\in\Gamma(X/Y)$ and $x^{\ast}\in X^{\ast}$ such that
$f=h\circ\pi+x^{\ast}$ if and only if $\operatorname*{qri}\overline
{\operatorname*{dom}f^{\ast}}\neq\emptyset$. In such a case, $x^{\ast}%
\in\operatorname*{qri}\overline{\operatorname*{dom}f^{\ast}}$ and $Y=L_{f}$.

\emph{(iii)} Assume that $(X,\left\langle \cdot,\cdot\right\rangle )$ is a
Hilbert space and $\operatorname*{qri}\overline{\operatorname*{dom}f^{\ast}%
}\neq\emptyset$. Then there exist a unique closed linear subspace $Y$ of $X$,
a unique essentially directionally coercive function $c\in\Gamma(Z)$ with
$Z:=Y^{\perp}$, and a unique $v\in Y$ such that $f=c\circ\Pr_{Z}+\left\langle
\cdot,v\right\rangle $, where $\Pr_{Z}$ is the orthogonal projection of $X$
onto $Z$. More precisely, $Y=L_{f}$, $c=h|_{Z}$ and $v:=\Pr_{Y}(x^{\ast})$ for
some (any) $x^{\ast}\in\overline{\operatorname*{dom}f^{\ast}}$, where $\Pr
_{Y}=I-\Pr_{Z}$.
\end{theorem}

Proof. (i) Take $x^{\ast}\in X^{\ast}$ and consider $Y:=L_{x^{\ast}}:=\{u\in
X\mid f_{\infty}(\pm u)=\left\langle \pm u,x^{\ast}\right\rangle \}$. Then $Y$
is closed linear subspace of $X$ by Lemma \ref{lem-Fact 5}. Using Proposition
\ref{p-Fact 6}~(a) we get $h\in\Gamma(X/Y)$ such that $f=h\circ\pi+x^{\ast}$
and $h_{\infty}(\pm\widehat{u})=0~\Rightarrow$ $\widehat{u}=\widehat{0}$;
hence $h$ is not constant on any line by \ref{r-a4a}. The other conclusions
follow from Proposition \ref{p-Fact 6}~(d) and (e).

(ii) The assertion is a consequence of (i) and Proposition \ref{p-Fact 6}~(c).

(iii) We identify $X^{\ast}$ with $X$ by Riesz theorem; then, for $Y$ a closed
linear subspace of $X$, the natural projection $\pi$ of $X$ onto $X/Y$ becomes
the orthogonal projection of $X$ on $Y^{\perp}$.

Assuming that $f=c\circ\Pr_{Z}+\left\langle \cdot,v\right\rangle $ with
$c\in\Gamma(Z)$ essentially directionally coercive and $v\in Y$ (less the
uniqueness), then $c=h+\left\langle \cdot,z\right\rangle $ with $h\in
\Gamma(Z)$ directionally coercive ($\Leftrightarrow\lbrack h_{\infty}%
\leq0]=\{0\}$) and $z\in Z$, whence
\[
f=h\circ\Pr\nolimits_{Z}+\left\langle \cdot,z+v\right\rangle .
\]
Having in view Proposition \ref{p-Fact 4}, because $\Pr_{Z}$ is onto, one must
have $x^{\ast}:=\left\langle \cdot,z+v\right\rangle \in\operatorname*{qri}%
\overline{\operatorname*{dom}f^{\ast}}$ and $(Z^{\perp}=)$ $\ker\Pr
_{Z}=[f_{\infty}\leq x^{\ast}]$. Using Proposition \ref{p-Fact 2a}, one must
have $(Y=)$ $Z^{\perp}=L_{x^{\ast}}=L_{f}$, whence $Z=L_{f}^{\perp}%
=\overline{\operatorname*{span}}\left(  \overline{\operatorname*{dom}f^{\ast}%
}-\overline{\operatorname*{dom}f^{\ast}}\right)  $ by Proposition
\ref{p-Fact 2}~(b); in particular, we got the uniqueness of $Y$. In order to
get the uniqueness of $v$, let us consider $x_{1}^{\ast},x_{2}^{\ast}%
\in\overline{\operatorname*{dom}f^{\ast}}$. Then $x_{i}^{\ast}=u_{i}+v_{i}$
with $u_{i}\in Z$ and $v_{i}\in Y$ for $i=1,2$. It follows that $Z\ni
x_{1}^{\ast}-x_{2}^{\ast}=(u_{1}-u_{2})+(v_{1}-v_{2})$. Because $Z\cap
Y=\{0\}$, we obtain that $v_{1}=v_{2}$. This shows that $\Pr_{Y}%
(\overline{\operatorname*{dom}f^{\ast}})$ is a singleton $\{v\}$. Because
$v\in Y$ and $c\circ\Pr_{Z}=f-\left\langle \cdot,v\right\rangle $, we have
that $c(z)=c(\Pr_{Z}(z))=f(z)-\left\langle z,v\right\rangle =f(z)$ for $z\in
Z$, that is $c=f|_{Z}$. This proves the uniqueness of $c$ in the
representation $f=c\circ\Pr_{Z}+\left\langle \cdot,v\right\rangle $ with the
desired properties.

In what concerns the existence of $Y$, $c$ and $v$ with the desired
properties, we proceed as follows: Consider $x^{\ast}\in\operatorname*{qri}%
\overline{\operatorname*{dom}f^{\ast}}$ and $Y=L_{f}$ $(=L_{x^{\ast}})$; set
$Z:=Y^{\perp}$ $(=X/Y)$. By (ii) there exist $h\in\Gamma(Z)$ directionally
coercive and $x^{\ast}\in X^{\ast}$ $(=X)$ such that $f=h\circ\pi+x^{\ast}$
$(=h\circ\Pr_{Z}+x^{\ast})$. Take $c:=f|_{Z}$, $v:=\Pr_{Y}(x^{\ast})\in Y$ and
$z:=x^{\ast}-v\in$ $Z$. Then
\[
c(z^{\prime})=f(z^{\prime})=(h\circ\Pr\nolimits_{Z})(z^{\prime})+\langle
z^{\prime},x^{\ast}\rangle=h(z^{\prime})+\langle z^{\prime},z+v\rangle
=h(z^{\prime})+\langle z^{\prime},z\rangle\quad\forall z^{\prime}\in Z,
\]
that is $c=h+\left\langle \cdot,z\right\rangle $. Hence $c$ is essentially
directionally coercive and $f=c\circ\Pr_{Z}+\left\langle \cdot,v\right\rangle
$. \hfill$\square$

\medskip

Let us see the relationships among the sets and functions considered in
\cite[Ths.\ 4--6]{Aza:19} and those introduced and used previously. As in
\cite{Aza:19}, in the sequel $X$ $(:=Z)$ is a Banach space.

As already observed at the beginning of this section, for $f\in\Gamma(X)$ one
has $\operatorname{Im}\partial f\subset\operatorname*{dom}f^{\ast}%
\subset\overline{\operatorname*{dom}f^{\ast}}$ and $(L_{f})^{\perp}%
=\overline{\operatorname*{span}}\left(  \operatorname{Im}\partial
f-\operatorname{Im}\partial f\right)  $ in the present framework.

\smallskip We begin with \cite[Th.\ 5]{Aza:19}. We observed before Lemma
\ref{lem-Fact a} that $Y_{f}=\{v\in Z\mid f_{\infty}(\pm v)=\left\langle \pm
v,\xi_{0}\right\rangle \}$ for $\xi_{0}\in\partial f(z_{0})\subset
\operatorname*{dom}f^{\ast}\subset\partial f_{\infty}(0)$, and so
$Y_{f}=L_{\xi_{0}}=L_{f}$ by (\ref{r-a15}).

Applying Theorem \ref{t-aza}~(i) for $x^{\ast}\in\operatorname{Im}\partial f$
one obtains a weaker version of \cite[Th.\ 5]{Aza:19}; more precisely $c_{f}$
attains its infimum on $Z/Y_{f}$ instead of taking values in $\mathbb{R}_{+}$.
In fact, in general it is not possible to obtain $c_{f}:Z/Y_{f}\rightarrow
\mathbb{R}_{+}$ with the desired properties in \cite[Th.\ 5]{Aza:19}. Indeed,
take $f(x):=\left\Vert x\right\Vert -1$; hence $f^{\ast}=\iota_{U_{X^{\ast}}%
}+1$. In this case $L_{f}=\{0\}$ and $\inf h_{x^{\ast}}=-f^{\ast}(x^{\ast
})=-1$ for every $x^{\ast}\in\operatorname*{dom}f^{\ast}=U_{X^{\ast}}$.

\medskip

Theorem 6 from \cite{Aza:19} is a obtained from Theorem \ref{t-aza}~(ii)
taking $\ell:=x^{\ast}\in\operatorname*{qri}\overline{\operatorname*{dom}%
f^{\ast}}$; such an $x^{\ast}$ exists by Proposition \ref{p-Fact 2c} because
$Z$ is separable. In fact $c_{f}$ is even directionally coercive with this
choice of $\ell$. \medskip

With respect to \cite[Th.\ 4]{Aza:19}, first observe that the space $X_{f}$
$[:=\operatorname*{span}(\operatorname{Im}\partial f-\operatorname{Im}\partial
f)]$ is not closed in general; for example, for the function $\varphi:\ell
_{2}\rightarrow\mathbb{R}$ defined by $\varphi(x):=\sum_{n=1}^{\infty
}\left\vert x_{n}\right\vert ^{2}/2^{n}$ (on page 2 of \cite{Aza:19}) one has
$\operatorname{Im}\partial\varphi=\{(y_{n})_{n\geq1}\in\ell_{2}\mid(2^{n}%
y_{n})_{n\geq1}\in\ell_{2}\}$ and $\operatorname*{dom}\varphi^{\ast}%
=\{(y_{n})_{n\geq1}\in\ell_{2}\mid(2^{n/2}y_{n})_{n\geq1}\in\ell_{2}\}$. Hence
$X_{f}$ is not closed for $f:=\varphi$. Replacing $X_{f}$ by its closure in
\cite[Th.\ 4]{Aza:19}, its conclusion follows from Theorem \ref{t-aza}~(iii)
because $\operatorname*{qri}\overline{\operatorname*{dom}f^{\ast}}$ is
nonempty, $Z$ being separable.

\medskip We end this note with an example which could be useful for providing
(counter-) examples.

\begin{example}
\label{ex-Fact 6b}Let $X$ be a normed vector space and $C\subset X^{\ast}$ be
a nonempty $w^{\ast}$-closed convex set. Then $\varphi_{C}:=(\tfrac{1}%
{2}\left\Vert \cdot\right\Vert ^{2})\nabla s_{C}$ with $s_{C}(x):=\sup
_{x^{\ast}\in C}\left\langle x,x^{\ast}\right\rangle $ for $x\in X$ is a
real-valued continuous convex function such that $\operatorname*{dom}%
\varphi_{C}^{\ast}=C$ and $(\varphi_{C})_{\infty}=s_{C}$. Here $h_{1}\nabla
h_{2}$ denotes the convolution of the functions $h_{1},h_{2}:X\rightarrow
\overline{\mathbb{R}}$ and is defined by $(h_{1}\nabla h_{2})(x):=\inf
\{h_{1}(x_{1})+h_{2}(x_{2})\mid x_{1},x_{2}\in X$, $x_{1}+x_{2}=x\}$.
\end{example}

Proof. Clearly, $s_{C}$ is a proper sublinear lsc function with $\psi^{\ast
}=\iota_{C}$. By \cite[Exer.\ 3.11~1)]{Zal:02} we have that $\varphi_{C}$ is a
continuous convex function such that $\varphi_{C}\leq\tfrac{1}{2}\left\Vert
\cdot\right\Vert ^{2}$, while from \cite[Th.\ 2.3.1~(ix)]{Zal:02},
$\varphi_{C}^{\ast}=(\tfrac{1}{2}\left\Vert \cdot\right\Vert ^{2})^{\ast
}+s_{C}^{\ast}=\tfrac{1}{2}\left\Vert \cdot\right\Vert ^{2}+\iota_{C}$. Hence
$\operatorname*{dom}\varphi_{C}^{\ast}=C$, whence $(\varphi_{C})_{\infty
}=s_{C}$ by (\ref{r-a5}). \hfill$\square$

\medskip Notice that taking $X:=\ell_{2}(\Gamma)$ and $C:=X_{+}$ as defined in
Remark \ref{rem2}, and $f$ the function defined in \cite[Ex.\ 7]{Aza:19}, then
$f=2\varphi_{C}$, where $\varphi_{C}$ is defined in Example \ref{ex-Fact 6b}.
Then $\operatorname*{dom}f^{\ast}=X_{+}$. So $L_{f}=(X_{+}-X_{+})^{\perp
}=\{0\}$ which shows that $f$ is not constant on any line [by Theorem
\ref{t-aza}~(i)]; moreover, if $\Gamma$ is uncountable, then
$\operatorname*{qri}\overline{\operatorname*{dom}f^{\ast}}=\operatorname*{qi}%
\overline{\operatorname*{dom}f^{\ast}}=\emptyset$ by Remark \ref{rem2}, and so
$f$ is not essentially directionally coercive by Theorem \ref{t-aza}~(ii). So,
the conclusions of \cite[Ex.\ 7]{Aza:19} are confirmed.

\end{document}